\documentstyle[amscd,amssymb,verbatim,epsf]{amsart}

\begin{document}
\theoremstyle{plain}
\newtheorem{Thm}{Theorem}
\newtheorem{Cor}{Corollary}
\newtheorem{Ex}{Example}
\newtheorem{Con}{Conjecture}
\newtheorem{Main}{Main Theorem}
\newtheorem{Lem}{Lemma}
\newtheorem{Prop}{Proposition}

\theoremstyle{definition}
\newtheorem{Def}{Definition}
\newtheorem{Note}{Note}

\theoremstyle{remark}
\newtheorem{notation}{Notation}
\renewcommand{\thenotation}{}

\errorcontextlines=0
\numberwithin{equation}{section}
\renewcommand{\rm}{\normalshape}%

\title[Surfaces of Constant Width]%
   {On C$^2$-smooth Surfaces of Constant Width}

\author{Brendan Guilfoyle}
\address{Brendan Guilfoyle\\
          Department of Mathematics and Computing \\
          Institute of Technology, Tralee \\
          Clash \\
          Tralee  \\
          Co. Kerry \\
          Ireland.}
\email{brendan.guilfoyle@@ittralee.ie}
\author{Wilhelm Klingenberg}
\address{Wilhelm Klingenberg\\
 Department of Mathematical Sciences\\
 University of Durham\\
 Durham DH1 3LE\\
 United Kingdom.}
\email{wilhelm.klingenberg@@durham.ac.uk }

\keywords{convex geometry, constant width, line congruence}
\date{April 24th, 2007}

\begin{abstract}
A number of results for C$^2$-smooth surfaces of constant width in Euclidean 3-space ${\mathbb{E}}^3$ are obtained. In particular, an 
integral inequality for constant width surfaces is established. This is used to prove that the ratio of volume to 
cubed width of a constant width surface is reduced by shrinking it along its normal lines. 
We also give a characterization of surfaces of constant width that have rational support function.  

Our techniques, which are complex differential geometric in nature, allow us to construct explicit smooth 
surfaces of constant width in ${\mathbb{E}}^3$, and their focal sets. They also allow for easy construction of tetrahedrally symmetric surfaces of constant width.
\end{abstract}

\maketitle

\section{Introduction}

The width of a closed convex subset of Euclidean ${\mathbb{E}}^n$ is  the
distance between parallel supporting planes, which is a map w: S$^{n-1}\rightarrow{\mathbb{R}}$. 
Subsets of constant width have been the studied in the 
context of convex geometry for many decades - see \cite{chakgrom} and references therein.

The purpose of this note is to bring some new differential geometric tools to bear on the construction of
subsets of constant width in ${\mathbb{E}}^3$, which we identify with their boundary surface. 
The nature of these tools are such that this boundary will be at least C$^2$-smooth.  

Our interest in developing these tools is two-fold. On the one hand, the Blaschke-Lebesgue problem of 
finding the convex body of 
fixed constant width of minimal volume in ${\mathbb{E}}^n$ remains open in dimensions greater than 2. While such a 
minimizer is not likely to be C$^2$-smooth, let alone smooth, it should be possible to approximate the minimizer by a 
constant width surface 
with degree k rational support function and induct on k. On the other hand, bodies of constant width play a central
role in  research on the potential theory of the farthest point distance function. Indeed,
a conjecture of Pritsker is complimentary to the Blaschke-Lebesgue problem in dimension 2 and open in 
higher dimensions \cite{gardnet} \cite{laugprits}.

Firstly, we establish an integral inequality for C$^2$-smooth surfaces of constant width (Theorem \ref{t:ineq}). 
If we move a  surface of 
constant width a fixed distance along its normal lines, the resulting ``parallel'' surface also has constant width.
The integral involved is invariant under such a shift and it is really from this perspective that our geometric
approach arises.  

We utilise the inequality to prove that, given a surface of constant width, shrinking
the surface along its inward pointing normal line reduces the volume with respect to its cubed 
width (Theorem \ref{t:flow}). 
Thus if we seek to solve the Blaschke-Lebesgue problem within a family of parallel constant width surfaces, 
we must squeeze the surface down along its normal as far as possible. The obstruction here is loss of convexity of 
the surface, which can also be characterized as the point at which the surface first touches its focal set.
Our techniques also allow for the computation of focal sets of arbitrary line congruences \cite{gak3}, which we can
then utilise.

Secondly, we characterize surfaces of constant width with rational support function. In particular, we prove that
the denominator must
satisfy a generalised palindromic condition utilising the antipodal map on S$^2$. 
Working within the rational support function class, 
we find evidence that the minimal volume obtained by shrinking along the normal is independent of the numerator of the 
support function. 

Finally, it is a conjecture of Danzer \cite{grubsch} that the minimizer of the Blaschke-Lebesgue problem in dimension 3
must have tetrahedral symmetry.
In fact, our techniques give a natural way to construct surfaces of constant width exhibiting any discrete symmetry: one
simply takes an arbitrary surface of constant width and sums over the elements of the group.
The result, which is also of constant width, has the symmetry, and in many cases, has smaller volume to width
ratio.

In the next section we summarise the pertinent geometric details culled from \cite{gak2} \cite{gak3}. In section 3 we 
apply this work to constant width surfaces, while the final section discusses examples of the construction in
detail.

\vspace{0.2in}

\section{Geometric Background}

\subsection{The Space of Oriented Lines}\label{s:1.1}
We start with 3-dimensional Euclidean space ${\mathbb{E}}^3$ and fix standard coordinates
($x^1,x^2,x^3$). In what follows we combine the first two coordinates
to form a single complex coordinate $z=x^1+ix^2$, set $t=x^3$ and refer to coordinates ($z,t$)
on ${\mathbb{E}}^3$.

\vspace{0.1in}
\setlength{\epsfxsize}{5.0in}
\begin{center}
   \mbox{\epsfbox{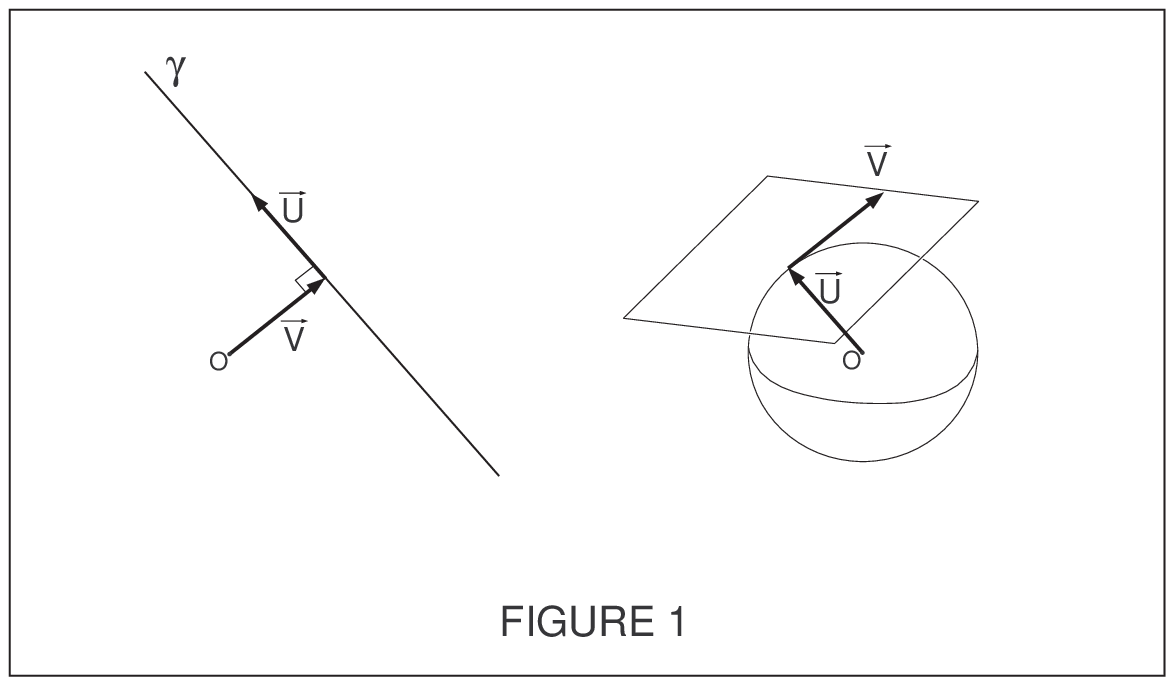}}
\end{center}
\vspace{0.1in}

Let ${\mathbb{L}}$ be the set of oriented lines, or {\it rays}, in Euclidean space
${\mathbb{E}}^3$. Such a line $\gamma$ is uniquely determined by its unit direction vector
$\vec{\mbox{U}}$ and the vector $\vec{\mbox{V}}$ joining the origin to the point on the line
that lies closest to the origin. That is, 

\[
\gamma=\{\;\vec{\mbox{V}}+r\vec{\mbox{U}}\in{\mathbb{E}}^3\;|\;r\in{\mathbb{R}}\;\},
\]
where $r$ is an affine parameter along the line.

By parallel translation, we move $\vec{\mbox{U}}$ to the origin 
and $\vec{\mbox{V}}$ to the head of $\vec{\mbox{U}}$. Thus, we obtain a vector that is tangent to
the unit 2-dimensional sphere in ${\mathbb{E}}^3$. 
The mapping is one-to-one and so it identifies the space of oriented lines 
with the tangent bundle of the 2-sphere ${\mbox{T\:S}}^2$ (see Figure 1).

\[
{\mathbb{L}}=\{\;(\vec{\mbox{U}},\vec{\mbox{V}})\in{\mathbb{E}}^3\times{\mathbb{E}}^3\;
               |\;\quad |\vec{\mbox{U}}|=1\quad\vec{\mbox{U}}\cdot\vec{\mbox{V}}=0\;\}.
\]

\vspace{0.2in}

\subsection{Coordinates on ${\mathbb{L}}$}\label{s:1.2}

The space ${\mathbb{L}}$ is a 4-dimensional manifold and the above identification gives a 
natural set of local complex coordinates.
Let $\xi$ be the local complex coordinate on the unit 2-sphere in ${\mathbb{E}}^3$
obtained by stereographic projection from the south pole. 

In terms of the standard spherical polar angles $(\theta,\phi)$, we have
 $\xi=\tan(\frac{\theta}{2})e^{i\phi}$. We convert from coordinates ($\xi,\bar{\xi}$) back to
($\theta,\phi$) using
\[
\cos\theta={\textstyle{\frac{1-\xi\bar{\xi}}{1+\xi\bar{\xi}}}}
\qquad
\sin\theta={\textstyle{\frac{2\sqrt{\xi\bar{\xi}}}{1+\xi\bar{\xi}}}}
\qquad
\cos\phi={\textstyle{\frac{\xi+\bar{\xi}}{2\sqrt{\xi\bar{\xi}}}}}
\qquad
\sin\phi={\textstyle{\frac{\xi-\bar{\xi}}{2i\sqrt{\xi\bar{\xi}}}}}.
\]

This can be extended to complex
coordinates $(\xi,\eta)$ on ${\mathbb{L}}$ minus the tangent space over the south
pole, as follows. First note that a tangent vector $\vec{\mbox{X}}$ to the 2-sphere can  
always be expressed as a linear combination of the tangent vectors generated by $\theta$ and $\phi$:
\[
\vec{\mbox{X}}=X^\theta\frac{\partial}{\partial\theta}+X^\phi\frac{\partial}{\partial\phi}.
\]
In our complex formalism, we have the natural complex tangent vector
\[
\frac{\partial}{\partial\xi}=\cos^2({\textstyle{\frac{\theta}{2}}})
    \left(\frac{\partial}{\partial\theta}
       -\frac{i}{2\cos({\textstyle{\frac{\theta}{2}}})\sin({\textstyle{\frac{\theta}{2}}})}\frac{\partial}{\partial\phi}
    \right) e^{-i\phi},
\]
and any real tangent vector can be written as
\[
\vec{\mbox{X}}=\eta\frac{\partial}{\partial\xi}+\bar{\eta}\frac{\partial}{\partial\bar{\xi}},
\]
for a complex number $\eta$. We identify the real tangent vector $\vec{\mbox{X}}$ on the 2-sphere 
(and hence the ray in ${\mathbb{E}}^3$) with the two complex numbers ($\xi,\eta$). Loosely speaking,
$\xi$ determines the direction of the ray, and $\eta$ determines its perpendicular distance
vector to the origin - complex representations of the vectors $\vec{\mbox{U}}$ and $\vec{\mbox{V}}$.

The coordinates ($\xi,\eta$) do not cover all of ${\mathbb{L}}$ - they omit all of the lines 
pointing directly downwards. However, the construction can also be carried out 
using stereographic projection from the north pole, yielding a coordinate system that 
covers all of ${\mathbb{L}}$ except for the lines pointing directly upwards. Between these two 
coordinate patches the whole of the space of oriented lines is covered. In what follows we work in the patch 
that omits the south direction.

\vspace{0.2in}

\subsection{The Correspondence Space}\label{s:1.4}

Geometric data will be transferred between ${\mathbb{E}}^3$ and ${\mathbb{L}}$ by use of a correspondence space.

\begin{Def}
The map $\Phi:{\mathbb{L}}\times{\mathbb{R}}\rightarrow{\mathbb{E}}^3$ is defined to take 
$((\xi,\eta),r)\in{\mathbb{L}}\times{\mathbb{R}}$ to the point  in ${\mathbb{E}}^3$ on the oriented line ($\xi,\eta$) 
that lies a distance $r$ from the point on the line closest to the origin (see the right of Figure 2).
\end{Def}

The double fibration on the left gives us the correspondence between the points in ${\mathbb{L}}$ and oriented lines in ${\mathbb{E}}^3$:
we identify a point ($\xi,\eta$) in ${\mathbb{L}}$ with $\Phi\circ\pi_1^{-1}(\xi,\eta)\subset{\mathbb{E}}^3$, which is an oriented line.
Similarly, a point p in ${\mathbb{E}}^3$ is identified with the 2-sphere $\pi_1\circ\Phi^{-1}(p)\subset{\mathbb{L}}$, which consists of all of the 
oriented lines through the point p.

\vspace{1.0in}
\begin{center}

\unitlength0.5cm

\begin{picture}(16,6)
\put(0.1,5.8){$\pi_1$}
\put(0.4,7.5){${\mathbb{L}}\times{\mathbb{R}}$}
\put(2,7){\vector(1,-1){2.8}}
\put(3.5,5.8){$\Phi$}
\put(0.8,3){${\mathbb{L}}$}
\put(1,7){\vector(0,-1){3}}
\put(5,3.4){${\mathbb{E}}^3$}
\put(7,-1){FIGURE 2}
\put(10,0){\mbox{\epsfbox{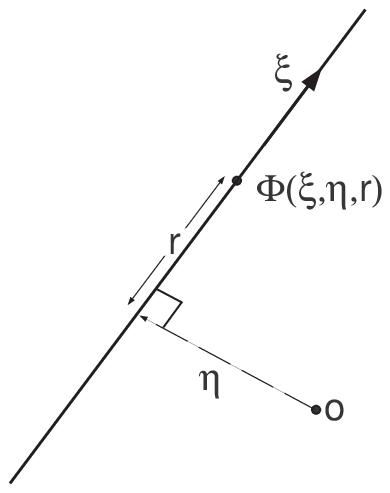}}}
\end{picture}

\end{center}
\vspace{0.2in}

The map $\Phi$ is of crucial importance when describing surfaces in ${\mathbb{E}}^3$ and
has the following coordinate expression:

\begin{Prop}\cite{gak2}
If $\Phi(\xi,\eta,r)=(z(\xi,\eta,r),t(\xi,\eta,r))$, then:
\begin{equation}\label{e:coord}
z=\frac{2(\eta-\overline{\eta}\xi^2)+2\xi(1+\xi\overline{\xi})r}{(1+\xi\overline{\xi})^2}
\qquad\qquad
t=\frac{-2(\eta\overline{\xi}+\overline{\eta}\xi)+(1-\xi^2\overline{\xi}^2)r}{(1+\xi\overline{\xi})^2},
\end{equation}
where $z=x^1+ix^2$, $t=x^3$ and ($x^1$, $x^2$, $x^3$) are Euclidean
coordinates in ${\mathbb{E}}^3$. 
\end{Prop}

\vspace{0.2in}

\subsection{Line Congruences}\label{s:3.1}

\begin{Def}
A {\it line congruence} is a 2-parameter family of oriented lines in ${\mathbb{E}}^3$. 
\end{Def}

From our perspective a line congruence is a surface $\Sigma$ in ${\mathbb{L}}$. In practice,
this will be given locally by a map 
${\mathbb{C}}\rightarrow{\mathbb{L}}:\mu\mapsto(\xi(\mu,\bar{\mu}),\eta(\mu,\bar{\mu}))$. A convenient 
choice of parameterization will depend upon the situation. In our case, the line congruences can be
parameterised by their directions. Thus we have $\xi\rightarrow(\xi,\eta=F(\xi,\bar{\xi}))$ and we label the
following combination of slopes
\begin{equation} \label{d:slopes}
\psi=(1+\xi\bar{\xi})^2\frac{\partial}{\partial \xi}\left(\frac{F}{(1+\xi\bar{\xi})^2}\right)
\qquad\qquad
\sigma=-\frac{\partial \bar{F}}{\partial \xi}.
\end{equation}

Given a line congruence $\Sigma\subset{\mathbb{L}}$, a map
$r:\Sigma\rightarrow {\mathbb{R}}$ determines a map
$\Sigma\rightarrow{\mathbb{E}}^3$ by
$(\xi,\eta)\mapsto\Phi((\xi,\eta),r(\xi,\eta))$ for $(\xi,\eta)\in\Sigma$. In other words,
we pick out one point on each line in the congruence (see Figure 3).

\vspace{0.1in}
\setlength{\epsfxsize}{5.0in}
\begin{center}
   \mbox{\epsfbox{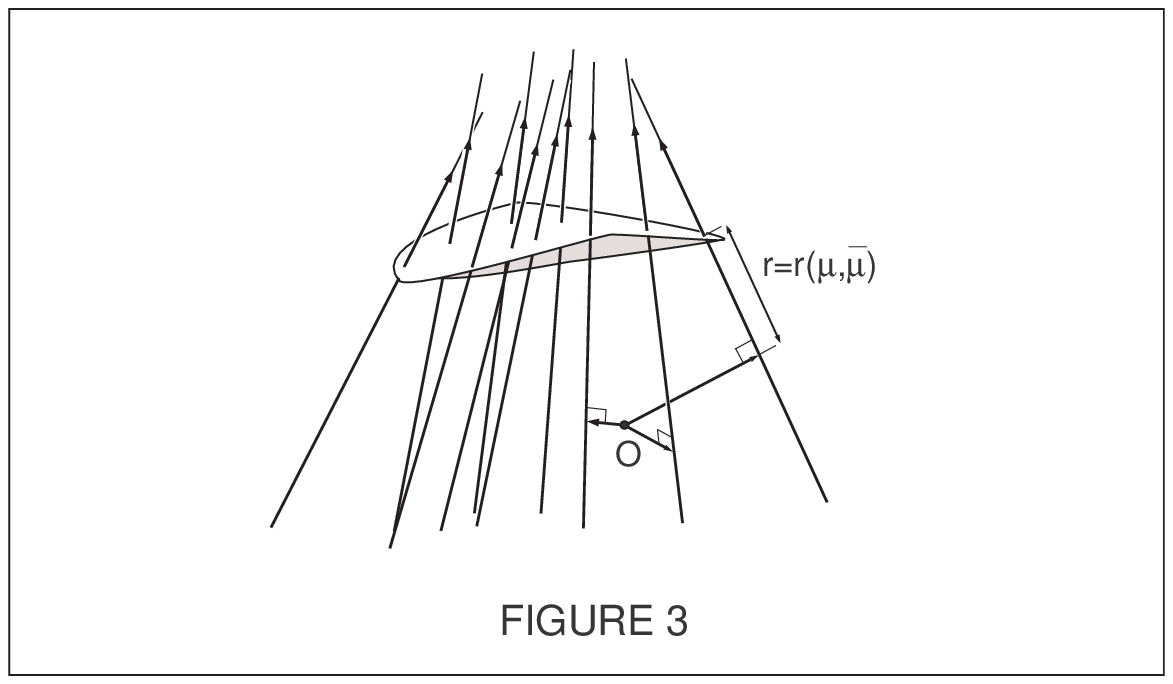}}
\end{center}
\vspace{0.1in}

For this surface to be 
orthogonal to the lines in ${\mathbb{E}}^3$, the complex function $F$ must
satisfy a certain condition:
  
\begin{Thm}\cite{gak2}
A line congruence $(\xi,\eta=F(\xi,\bar{\xi}))$ is orthogonal to a surface in ${\mathbb{E}}^3$ iff there exists a
real function $r(\xi,\bar{\xi})$ satisfying:
\begin{equation}\label{e:intsur}
\frac{\partial r}{\partial \bar{\xi}}=\frac{2F}{(1+\xi\bar{\xi})^2}.
\end{equation}

If there exists one solution, there exists a 1-parameter family generated by a real constant of integration.
The function r is the distance from the surface to the point on the normal line closest to the origin. 
\end{Thm}

\vspace{0.1in}

The surface can be reconstructed in ${\mathbb{E}}^3$ from this data be inserting $r=r(\xi,\bar{\xi})$ and 
$\eta=F(\xi,\bar{\xi})$ in equations (\ref{e:coord}). Note that condition (\ref{e:intsur}) implies that
the slope $\psi$ in (\ref{d:slopes}) is real.

\vspace{0.2in}

\subsection{Focal Points of a Line Congruence}\label{s:6.1}

Suppose we have a line congruence $\Sigma$ parameterized by its direction $\xi\rightarrow(\xi,\eta=F(\xi,\bar{\xi}))$.

\begin{Def}
A point $p\in{\mathbb{E}}^3$ on a line $\gamma$ in the line congruence $\Sigma$ is a {\it focal point} if the jacobian 
of the transformation $(\xi,r)\rightarrow\Phi((\xi,F(\xi,\bar{\xi})),r)$ vanishes at $p$.

The set of focal points of a line congruence $\Sigma$ generically form surfaces
in ${\mathbb{E}}^3$, which are referred to as the {\it focal surfaces} of $\Sigma$. 
\end{Def}

\begin{Thm}\label{t:focs}\cite{gak3}
The focal set of the parametric line congruence $\Sigma$ which is normal to a closed convex surface is given by
\[
r=r_{\pm}(\xi,\bar{\xi})=-\psi\pm |\sigma|,
\]
where the slopes $\psi$ and $\sigma$ are given by equation (\ref{d:slopes}). Thus on each there is either one or two
focal points. 
\end{Thm}

\vspace{0.2in}

\section{Surfaces of Constant Width}

\subsection{Oriented Normal lines}

Consider a closed convex body B in ${\mathbb{E}}^3$ with smooth boundary surface S. 
The set of oriented normal lines to S forms a line congruence that can be parameterized by the direction of
the normal. Thus the normals are given by a map $\xi\rightarrow(\xi,\eta=F(\xi,\bar{\xi}))$, and there
exists a real function $r(\xi,\bar{\xi})$ satisfying equation (\ref{e:intsur}).

\begin{Def}
The map $r:S^2\rightarrow{\mathbb{R}}$ is the distance of the tangent planes of S to the origin and is called the 
{\it support function} of S. If $\tau:S^2\rightarrow S^2$ is the antipodal map, the {\it width} of S is a 
function w:S$^2\rightarrow{\mathbb{R}}$ defined by $w=r+r\circ\tau$.
\end{Def}

\begin{Prop}
The oriented normals to a surface of constant width $w$ are given by 
$\xi\rightarrow(\xi,\eta=F(\xi,\bar{\xi}))$ where the lines have the reflection symmetry:
\[
F(\tau(\xi),\tau(\bar{\xi}))=-\frac{1}{\bar{\xi}^2}\overline{F(\xi,\bar{\xi})}.
\]
\end{Prop}
\begin{pf}
This follows from the fact that the antipodal map is $\tau(\xi)=-\bar{\xi}^{-1}$, differentiation
of the constant width condition and equation (\ref{e:intsur}).
\end{pf}

\vspace{0.2in}

\subsection{The Blaschke-Lebesgue Problem}

We now consider the volume of a closed convex body B in ${\mathbb{E}}^n$ with smooth boundary S. For ease of notation we 
denote the volume of B by ${\mbox {Vol}}(S)$, meaning, of course, the volume enclosed by S. Let $S^n_w$ be the 
round n-sphere of width w.

\begin{Def}
For a closed convex body in ${\mathbb{E}}^n$ of constant width w with boundary S, we define
\[
{\cal{I}}(S)=\frac{{\mbox {Vol}}(S)}{{\mbox {Vol}}(S^{n-1}_w)}.
\] 
\end{Def}

As a consequence of a well-known theorem of Bieberbach, the sphere S$^{n-1}$ maximises ${\cal{I}}$ in 
Euclidean ${\mathbb{E}}^n$. 
The problem of minimizing ${\cal{I}}$ was solved for $n=2$ by Blaschke and Lebesgue and turns out to be minimized by the
Reuleaux triangle \cite{bonnfen}.
While a number of shorter proofs have since been given for this result, the problem remains open for $n>2$.

For n$=3$, the smallest known example is a body with 
${\cal{I}}(S)=4-\frac{3\sqrt{3}}{2}\cos^{-1}\left(\frac{1}{3}\right)=0.801873619$ \cite{bonnfen}. On the other hand the 
best lower bound for ${\cal{I}}$ is $2(3\sqrt{6}-7)=0.696938456$ \cite{chak}, so a large gap remains. From here
on we consider only the case n$=3$.

In this context, a useful formula of Blaschke says that the volume enclosed by a surface S of constant width $w$ can be 
computed from the area $A(S)$ by
\begin{equation}\label{e:blaschke}
{\mbox {Vol}}(S)={\textstyle{\frac{1}{2}}}wA(S)-{\textstyle{\frac{1}{3}}}\pi w^3.
\end{equation}
Thus, to minimize the volume of the body we must minimize the surface area of the boundary. The following proposition 
gives an expression for the surface area in terms of the slopes of the normal line congruence:
 
\begin{Prop}
The surface area of a convex surface S with support function $r(\xi,\bar{\xi})$ is 
\begin{equation}\label{e:area}
A(S)=\int\int_{S^2}(r+\psi)^2-|\sigma|^2\frac{d\xi d\bar{\xi}}{(1+\xi\bar{\xi})^2},
\end{equation}
where, as before,
\[
\psi=(1+\xi\bar{\xi})^2\frac{\partial}{\partial \xi}\left(\frac{F}{(1+\xi\bar{\xi})^2}\right)
\qquad\qquad
\sigma=-\frac{\partial \bar{F}}{\partial \xi},
\]
and
\[
F={\textstyle{\frac{1}{2}}}(1+\xi\bar{\xi})^2\frac{\partial r}{\partial \bar{\xi}}.
\]
\end{Prop}
\begin{pf}
This follows immediately from the coordinate expression for a null basis found in the proof of Theorem 2 in \cite{gak2}.
\end{pf}

\vspace{0.1in}

We now prove an integral inequality for surfaces of constant width :

\begin{Thm}\label{t:ineq}
For a surface of constant width w with support function $r(\xi,\bar{\xi})$ 
\begin{equation}\label{e:ineq}
\int\int_{S^2}|\sigma|^2-(r-{\textstyle{\frac{1}{2}}}w+\psi)^2\frac{d\xi d\bar{\xi}}{(1+\xi\bar{\xi})^2}\geq0,
\end{equation}
where $\sigma$ and $\psi$ are given by (\ref{d:slopes}). Equality only occurs in the case of the 2-sphere of width w.
\end{Thm}
\begin{pf}
Given that $\tau(\xi)=-\bar{\xi}^{-1}$, a short computation shows that, for a surface of constant width w,
\[
r\circ\tau=w-r \qquad\qquad \psi\circ\tau=-\psi \qquad \qquad |\sigma\circ\tau|^2=|\sigma|^2.
\]
Now, since the area integral is invariant under the antipodal map we can average over the identity and the antipodal map
to get
\begin{align}
A(S)=&{\textstyle{\frac{1}{2}}}\int\int_{S^2}(r+\psi)^2-|\sigma|^2+(w-r-\psi)^2-|\sigma|^2
   \frac{d\xi d\bar{\xi}}{(1+\xi\bar{\xi})^2}\nonumber\\
&=\int\int_{S^2}(r-{\textstyle{\frac{1}{2}}}w+\psi)^2+{\textstyle{\frac{1}{4}}}w^2-|\sigma|^2
   \frac{d\xi d\bar{\xi}}{(1+\xi\bar{\xi})^2}\nonumber\\
&=\pi w^2-\int\int_{S^2}|\sigma|^2-(r-{\textstyle{\frac{1}{2}}}w+\psi)^2\frac{d\xi d\bar{\xi}}{(1+\xi\bar{\xi})^2}.
    \nonumber
\end{align}
By the theorem of Bieberbach mentioned earlier $A(S)\leq\pi w^2$ with equality iff S is the 2-sphere of width w. The
stated result follows from applying this to the above geometric identity.
\end{pf}

\vspace{0.1in}

We can apply this inequality as follows. If we move the points on a surface of constant width a fixed distance C along its
normal line we get another surface of constant width. Indeed, the support function changes by $r\rightarrow r+C$, the
width obviously changing by $w\rightarrow w+2C$. It is not immediately clear, however, how ${\cal{I}}$ changes under
such a shift. The following Theorem shows that it increases as C increases.

\begin{Thm}\label{t:flow}
Let $r=r_0$ be the support function of a C$^2$-smooth surface $S_0$ bounding a body of constant width w$_0$.
Let $S_C$ be the surface of constant width obtained from the support function $r=r_0+C$.
Then
\[
\frac{d }{dC}{\cal{I}}(S_C)\geq0.
\]
\end{Thm}
\begin{pf}
Since w$_0$ is the width of S$_0$, the width of S$_C$ is w$_0+2C$. We compute
\begin{align}
{\cal{I}}(S_C)=&\frac{{\mbox{Vol}}(S_C)}{{\mbox{Vol}}(S^2_{w_0+2C})} \nonumber \\
&=\left({\textstyle{\frac{1}{2}}}(w_0+2C)A(S_C)-{\textstyle{\frac{1}{3}}}\pi (w_0+2C)^3\right)
                \frac{6}{\pi(w_0+2C)^3}\nonumber \\
&=1-\frac{3}{\pi(w_0+2C)^2}\int\int_{S^2}|\sigma|^2-(r_0-{\textstyle{\frac{1}{2}}}w_0+\psi)^2\frac{d\xi d\bar{\xi}}{(1+\xi\bar{\xi})^2}\nonumber, 
\end{align}
where we have used Blaschke's formula (\ref{e:blaschke}) and the surface area formula (\ref{e:area}).

Now differentating we get
\[
\frac{d }{dC}{\cal{I}}(S_C)=\frac{6}{\pi(w_0+2C)}\int\int_{S^2}|\sigma|^2-(r_0-{\textstyle{\frac{1}{2}}}w_0+\psi)^2\frac{d\xi d\bar{\xi}}{(1+\xi\bar{\xi})^2}\geq0,
\]
as claimed.
\end{pf}

\vspace{0.1in}

Thus, to minimize ${\cal{I}}$ the constant width surface must be shrunk along its normal as far as possible, 
that is, until loss of convexity. Loss of convexity occurs when the surface comes into contact with its focal set
 \cite{gak3}. As we saw in the Theorem \ref{t:focs} this consists of two sets in ${\Bbb{E}}^3$ given by inserting
$r=-\psi\pm|\sigma|$ in (\ref{e:coord}). Thus, to minimize ${\cal{I}}$ we must find the minimum value for $C$ so that 
the surface just touches its focal set.

Focal sets are usually not smooth - they contain singular points which we refer to as cusps. At a point where the
focal set of a line congruence is smooth, the line is tangent to the focal set. Thus, it is clear that shrinking a 
convex surface S along its normal, the first point on the focal set that the surface S encounters will be a singular
point. In the next section we illustrate this.

\vspace{0.2in}

\subsection{Constant Width Surfaces with Rational Support}

\begin{Def}
A closed convex surface has {\it rational support} if the support function is of the form
\[
r=\frac{P(\xi,\bar{\xi})}{Q(\xi,\bar{\xi})}
\]
where $P$ and $Q$ are real-valued polynomials. Since P is real-valued, the degree of $\xi$ and $\bar{\xi}$ are equal, 
and we refer to this simply as the degree of P. Similarly, we have the degree of Q, and in order for the surface to
be closed we must have deg(P)$\leq$deg(Q). We also assume that $P\nmid Q$.
\end{Def}

We now characterize convex surface with rational support that are of constant width:

\begin{Thm}\label{t:ratsupp}
Consider a convex surface S with rational support, as above, with $deg(P)=n\leq deg(Q)=m$.  Then S
is of constant width w iff
\vspace{0.1in}
\begin{enumerate}
\item[(1)] 
\[
Q\left(-\frac{1}{\bar{\xi}},-\frac{1}{\xi}\right)=\frac{1}{K\;\xi^m\bar{\xi}^m}Q(\xi,\bar{\xi}),
\]
for some $K\in{\mathbb{R}}$.
\item[(2)] If 
\[
P(\xi,\bar{\xi})=\sum_{k,l=0}^mA_{kl}\xi^k\bar{\xi}^l
\qquad\qquad
Q(\xi,\bar{\xi})=\sum_{k,l=0}^mB_{kl}\xi^k\bar{\xi}^l,
\]
then
\[
A_{kl}+(-1)^{k+l}\;K\;A_{m-k\;m-l}=wB_{kl}.
\]
\end{enumerate}

\end{Thm}
\begin{pf}
We begin by complexifying 
\[
r(z_1,z_2)=\frac{P(z_1,z_2)}{Q(z_1,z_2)},
\]
for $z_1,z_2\in{\mathbb{C}}$. Thus P is of degree n in $z_1$ and $z_2$, while Q is of degree m in $z_1$ and $z_2$.
Define
\[
\tilde{P}(z_1,z_2)=z_1^nz_2^n\;P\left(-\frac{1}{z_2},-\frac{1}{z_1}\right)
\qquad
\tilde{Q}(z_1,z_2)=z_1^mz_2^m\;Q\left(-\frac{1}{z_2},-\frac{1}{z_1}\right).
\]
Now the antipodal map $\tau$ in holomorphic coordinates is $\tau(\xi)=-\bar{\xi}^{-1}$, so the constant 
width condition is
\[
\frac{P(z_1,z_2)}{Q(z_1,z_2)}+\frac{P(-z_2^{-1},-z_1^{-1})}{Q(-z_2^{-1},-z_1^{-1})}=w,
\]
or
\begin{equation}\label{e:cwcond1}
P(z_1,z_2)\tilde{Q}(z_1,z_2)+z_1^{m-n}z_2^{m-n}\tilde{P}(z_1,z_2)Q(z_1,z_2)=wQ(z_1,z_2)\tilde{Q}(z_1,z_2).
\end{equation}
Now for ($a_1,a_2$)$\in{\mathbb{C}}^2$ such that $Q(a,b)=0$ we have from the constant width condition (\ref{e:cwcond1}) 
that $P(a,b)\tilde{Q}(a,b)=0$. Since P and Q have no common factors, the complex curves in ${\mathbb{C}}^2$ given by
P$^{-1}$(0) and Q$^{-1}$(0) have no common components. Thus, except at a finite number of points,
\[
Q(a,b)=0 \qquad \Leftrightarrow \qquad \tilde{Q}(a,b)=0.
\]
But these are two polynomials of the same degree, and so we conclude that $Q(z_1,z_2)=K\;\tilde{Q}(z_1,z_2)$ for some
$K\in{\mathbb{C}}$. In fact, since the underlying polynomial is real-valued we see that $K\in{\mathbb{R}}$ and
\[
Q(z_1,z_2)=K\;\tilde{Q}(z_1,z_2)=K\;z_1^mz_2^mQ(-z_2^{-1},-z_1^{-1}).
\]
This establishes part (1).

To prove part (2) we compute
\begin{align}
w=&\frac{P(\xi,\bar{\xi})}{Q(\xi,\bar{\xi})}
+\frac{P\left(-\frac{1}{\bar{\xi}},-\frac{1}{\xi}\right)}{Q\left(-\frac{1}{\bar{\xi}},-\frac{1}{\xi}\right)}\nonumber\\
=&\left(\sum_{k,l=0}^mA_{kl}\xi^k\bar{\xi}^l+K\;\xi^m\bar{\xi}^m\sum_{k,l=0}^mA_{kl}(-\bar{\xi})^{-k}(-\xi)^{-l}\right)
     \left(\sum_{k,l=0}^mB_{kl}\xi^k\bar{\xi}^l\right)^{-1}\nonumber.
\end{align}
Thus,
\[
\sum_{k,l=0}^m\left(A_{kl}+(-1)^{k+l}K\;A_{n-k\;n-l}\right)\xi^k\bar{\xi}^l=w\sum_{k,l=0}^mB_{kl}\xi^k\bar{\xi}^l.
\]
Comparison of terms yields the result.
\end{pf}

\vspace{0.2in}

\section{Explicit Examples}

\subsection{Rotational Symmetry}

First consider the oriented normal lines to a convex surface that is rotationally symmetric about the $x^3$-axis.
It is not hard to see that the map $\xi\rightarrow(\xi,\eta=F(\xi,\bar{\xi}))$ determining this line congruence 
satisfies $F=G(R)e^{i\theta}$, where $G$ is a real function and $\xi=Re^{i\theta}$.

For rational support we have:

\begin{Cor}
Consider a convex surface S with rational support which is rotationally symmetric about the $x^3$-axis with
\[
P(R)=\sum_{k=0}^mA_{k}R^{2k}
\qquad\qquad
Q(R)=\sum_{k=0}^mB_{k}R^{2k}.
\]
Then S is of constant width w iff, after rescaling,
\begin{enumerate}
\item[(1)] Q is palindromic: $B_k=B_{m-k}$,
\item[(2)] P and Q satisfy
\[
A_{k}+A_{m-k}=wB_{k}.
\]
\end{enumerate} 
\end{Cor}

\vspace{0.1in}

We also have the following description of the focal sets:

\begin{Prop}\label{p:rotfoc}
The focal set of the oriented normals to a convex rotationally symmetric surface with support function 
$r=r(R)$ is given by the surface
\[
z={\textstyle{\frac{1}{2}}}\left(-R(1+R^2)\frac{d^2r}{dR^2}+(1-3R^2)\frac{dr}{dR}\right)e^{i\theta}
\]
\[
t={\textstyle{\frac{1}{4}}}\left(-(1-R^4)\frac{d^2r}{dR^2}-2R(3-R^2)\frac{dr}{dR}\right),
\]
and the line
\[
z=0
\qquad\qquad
t=-\frac{(1+R^2)^2}{4R}\frac{dr}{dR},
\]
where $z=x^1+ix^2$ and $t=x^3$, for standard coordinates $(x^1,x^2,x^3)$ on Euclidean 3-space.
\end{Prop}
\begin{pf}
This follows from Theorem \ref{t:focs} by imposing rotational symmetry and using
\[
\psi=r+\frac{(1+R^2)^2}{2R}\frac{d}{dR}\left(\frac{RG}{(1+R^2)^2}\right)
\qquad\qquad
\sigma=-\frac{1}{2}R\frac{d}{dR}\left(\frac{G}{R}\right)e^{-2i\theta},
\]
where
\[
G=\frac{1}{4}(1+R^2)^2\frac{dr}{dR}.
\]
\end{pf}

\vspace{0.1in}

Analogous results hold for focal sets of reflections off translation invariant surfaces \cite{gak4}.

The  singularities or cusps of the focal set of a rotationally symmetric surface are similarly described:

\begin{Prop}\label{p:rotcusp}
The cusps on the focal set of the oriented normals to a convex rotationally symmetric surface with support function 
$r=r(R)$ are solutions of the equation:
\begin{equation}\label{e:cusps}
(1+R^2)\frac{d^3r}{dR^3}+6\frac{d^2r}{dR^2}+6\frac{dr}{dR}=0.
\end{equation}
\end{Prop}
\begin{pf}
Cusps occur on the focal set given by the expressions in Proposition \ref{t:focs} when 
\[
\frac{d z}{dR}=0 \qquad\qquad {\mbox{and}} \qquad\qquad \frac{d t}{dR}=0.
\]
A straight-forward computation shows that these are equivalent to (\ref{e:cusps}).
\end{pf}

\vspace{0.2in}

\subsection{Example}\label{s:ex}

The support function
\[
r=\frac{a+bR^2+(3-b)R^4+(1-a)R^6}{(1+R^2)^3}+C,
\]
for a,b$\in{\mathbb{R}}$ gives a rotationally symmetric surface of constant width 1+2C. For $a=b-1$
this is a round sphere with centre ($0,0,b-{\textstyle{\frac{3}{2}}}$) and radius $C+{\textstyle{\frac{1}{2}}}$.

A straight-forward computation utilising (\ref{e:coord}) yields the parametric equation of the surface:
\[
x^1=\frac{\left[(a-b+2C+2)(3+R^4)R^2-(a-b-2C)(1+3R^4)\right]R\cos(\theta)}{(1+R^2)^4},
\]
\[
x^2=\frac{\left[(a-b+2C+2)(3+R^4)R^2-(a-b-2C)(1+3R^4)\right]R\sin(\theta)}{(1+R^2)^4},
\]
\[
x^3=\frac{\left[(a-C-1)R^8+(5a-b-2C-2)R^6+(6b-9)R^4+(5a-b+2C)R^2+a+C\right]}{(1+R^2)^4}.
\]

From our area formula (\ref{e:area}) we compute the volume and hence 
\[
{\cal{I}}=1-\frac{3(a-b+1)^2}{35(1+2C)^2}.
\]
Note again the sphere case when $a=b-1$.

We now compute the focal sets of the oriented normal lines, and Figure 4 illustrates the result. Since the surfaces
are all rotationally symmetric we only need consider a cross-section. The surface for different values of C
and the focal set, for a$=3$ and b$=3$ are shown. The focal set lying on the axis of symmetry is obtained from
$r=r_-=-\psi-|\sigma|$, while the triangular focal set is from $r=r_+=-\psi+|\sigma|$. We can see the loss of convexity 
once the surface crosses the cusps. Note that it hits all cusps at the same C-value.

To find these cusps we must solve equation (\ref{e:cusps}), which in our case works out to be 
\[
(a-b+1)R(R^2-3)(3R^2-1)=0.
\]
Since
$a-b+1\neq0$, we have cusps at $R=0$ and $R=\sqrt{3}$ and their antipodes. To find the $C$ at which the surface
just touches the cusps we compute
\[
r(0)-r_+(0)={\textstyle{\frac{1}{2}}}(-a+b+2C) 
\qquad\qquad
r(\sqrt{3})-r_+(\sqrt{3})={\textstyle{\frac{1}{2}}}(-a+b+2C) .
\]
The first point of contact with the focal set occurs when these vanish.
Thus the C value that minimizes ${\cal{I}}$ is $C=(a-b)/2$ and this value then works out to be 
${\cal{I}}=32/35=0.914285724$. 

\vspace{0.1in}
\setlength{\epsfxsize}{5.0in}
\begin{center}
   \mbox{\epsfbox{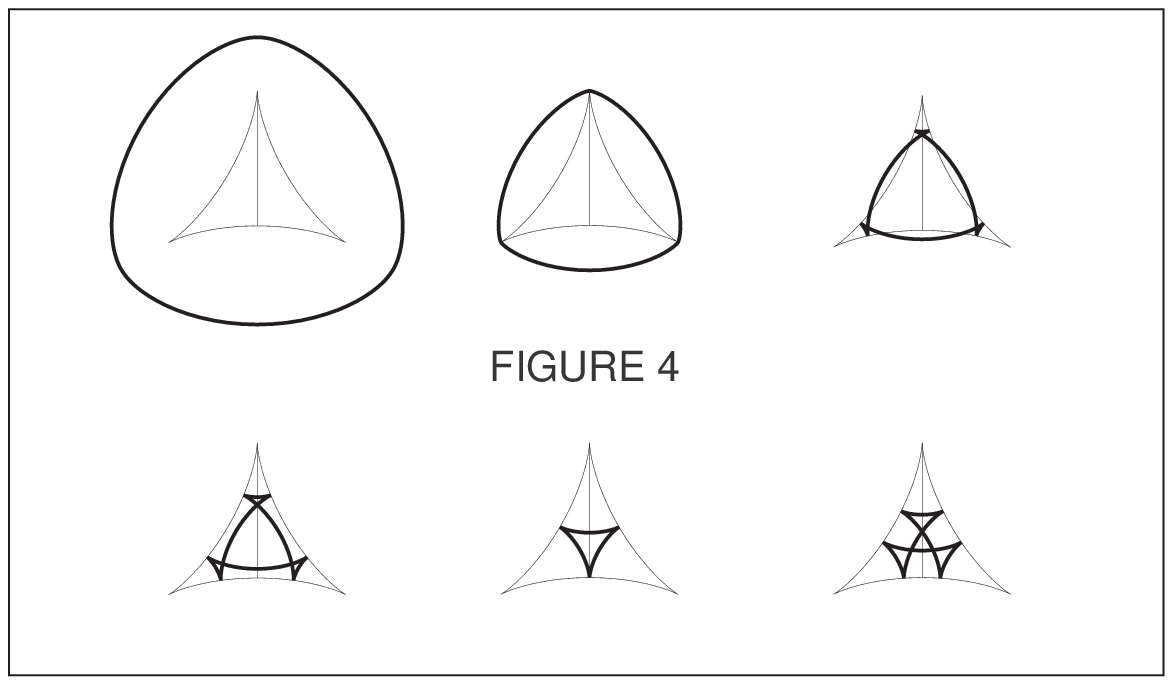}}
\end{center}
\vspace{0.1in}

It is remarkable that this value is independent of both a and b. We have a two-parameter
family of surfaces of constant width, but once they are shrunk along their normals they all yield surfaces 
enclosing the same volume. In fact, this property persists for higher powers of the denominator:

\begin{Prop}
Consider the constant width surfaces S given by
\[
r=\frac{P(|\xi|^2)}{(1+\xi\bar{\xi})^k},
\]
where the coefficients of P satisfy the conditions in Theorem \ref{t:ratsupp}. 

Then, ${\cal{I}}(S)=32/35$ for k=3,4,5,6,7.
\end{Prop}

\vspace{0.1in}

While induction on k in the above proposition is difficult to implement, we conjecture it should hold for all k. 
In fact, on the evidence of a large number of numerical experiments, we conjecture:

\vspace{0.1in}
\noindent{\bf Conjecture}:

Consider a constant width surfaces S with rational support function r. Then the functional ${\cal{I}}$ of the constant width surface 
obtained by shrinking the surface as far as possible along its normal lines is independent of the numerator of 
r.

\vspace{0.2in}

\subsection{Discrete Symmetries}

Consider a discrete subgroup of isometries ${\cal{G}}\subset$ O(3), and suppose that $r_0$ is the support function
of a surface of constant width w.

\begin{Prop}
The surface determined by the support function
\[
r(\xi,\bar{\xi})=\frac{1}{\#{\cal{G}}}\sum_{g\in {\cal{G}}}r_0(g(\xi),g(\bar{\xi})),
\]
is a surface of constant width w which is invariant under G.
\end{Prop}
\begin{pf}
This follows from the fact that the antipodal map commutes with elements of O(3).
\end{pf}

Applying this approach to the case of $r_0$ being equal to the support function in Example \ref{s:ex} and 
${\cal{G}}$ being the tetrahedral group, we can 
construct closed convex surfaces of constant width with tetrahedral symmetry. The results are shown in Figure 5,
where both a surface (left) and its focal set (right) is presented.

For this example, the minimum value of ${\cal{I}}$ obtained is approximately 0.8794644289, which is an improvement on the 
rotationally symmetric value.

\vspace{0.1in}
\setlength{\epsfxsize}{5.0in}
\begin{center}
   \mbox{\epsfbox{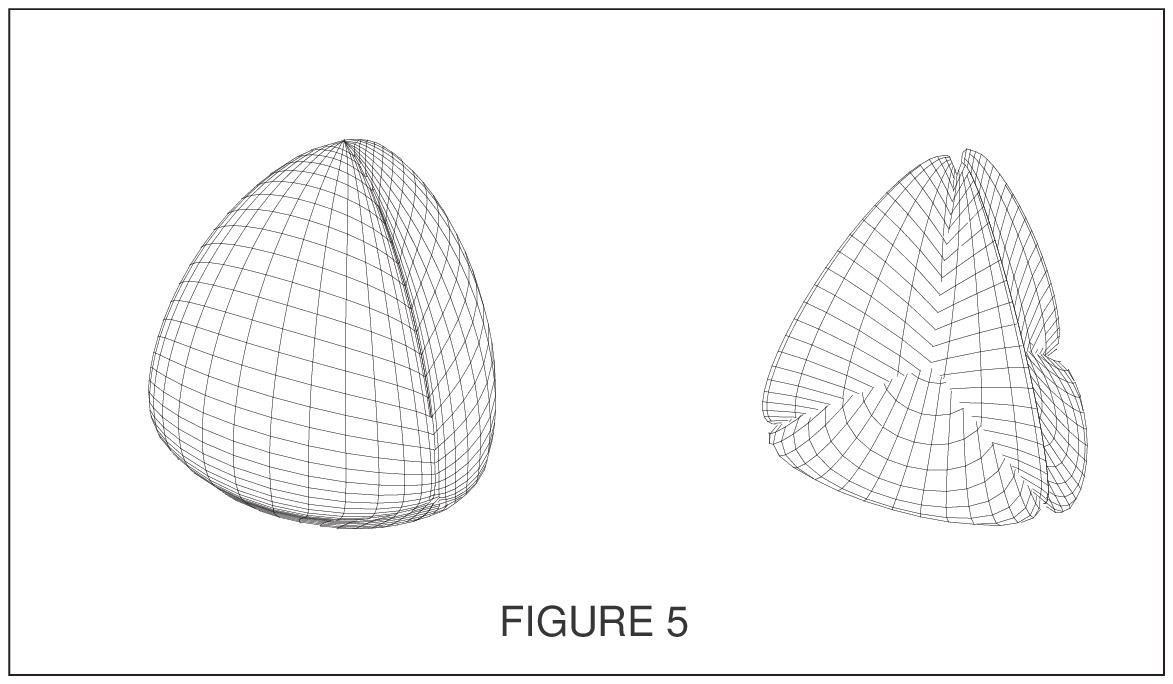}}
\end{center}
\vspace{0.1in}

\noindent{\bf Acknowledgement}: 

The authors would like thank Peter Giblin for bringing this topic to their attention. 
Part of this work was supported by the Research in Pairs Programme of the Mathematisches Forschungsinstitut Oberwolfach, 
Germany.

\end{document}